\documentclass[12pt]{amsart}
\usepackage{amsfonts}
\usepackage{amssymb}
\theoremstyle{definition}

\theoremstyle{remark}

\newcommand{\X}{\mathfrak X}

\newcommand{\ds}{\displaystyle}

\begin{document}

\centerline{\bf Comptes rendus de l'Academie bulgare des Sciences}

\centerline{\it Tome 36, No 10, 1983}

\vspace{0.6in}

\begin{flushright}
{\it MATH\'EMATIQUES
\\ G\'eometrie}
\end{flushright}

\vspace{0.2in}

\centerline{\large\bf FOUR-DIMENSIONAL CONFORMAL FLAT $QK_3$-MANIFOLDS}

\vspace{0.3in}
\centerline{\bf O. T. Kassabov}

\vspace{0.3in}
\centerline{\it (Submited by Academician H.Hristov on May 23, 1983)}

\vspace{0.2in}

Let $M$ be an almost Hermitian manifold with Riemannian metric $g$ and
almost complex structure $J$. We denote by $\nabla$ and $R$ the covariant
differentiation and the curvature tensor of $M$, respectively. Let 
$F(X,Y)=g(X,JY)$. Then $M$ is said to be a K\"ahler, or nearly K\"ahler, 
or almost K\"ahler, or quasi K\"ahler manifold, if $\nabla J=0$, or
$(\nabla_X J)X=0$, or $dF=0$, or 
$$
	(\nabla_X J)Y + (\nabla_{JX} J)JY =0 \ ,
$$
respectively. The corresponding classes of manifolds are denoted by $K$, $NK$,
$AK$, $QK$. For a given class $L$ of almost Hermitian manifolds the class $L_i$
is defined by the identity $(i)$, where
$$
	R(X,Y,Z,U)=R(X,Y,JZ,JU)  \leqno(1)
$$
$$
	R(X,Y,Z,U)=R(X,Y,JZ,JU)+ R(X,JY,Z,JU)+R(JX,Y,Z,JU)  \leqno(2)
$$
$$
	R(X,Y,Z,U)=R(JX,JY,JZ,JU) .  \leqno(3)
$$
It is known \cite{G}, that
$$
	\begin{array}{rcl}
		   &\subset \ NK \subset &  \\	 
		K  &                     &QK_2 \\
		   &\subset \ AK_2 \subset &
	\end{array} \qquad
	NK \cap AK = K
$$

A $2m$-dimensional almost Hermitian manifold ($m \ge 2$) is conformal flaf, iff
$$
	\begin{array}{r}
\ds R(X,Y,Z,U) = \frac{1}{2m-2}\{ g(X,U)S(Y,Z) - g(X,Z)S(Y,U)   \\
                                 +g(Y,Z)S(X,U) - g(Y,U)S(X,Z)\} \\
\ds -\frac{\tau}{(2m-1)(2m-2)} \{ g(X,U)g(Y,Z) - g(X,Z)g(Y,U) \}
	\end{array}    \leqno (1)
$$        
where $S$ (resp. $\tau$) is the Ricci tensor (resp. the scalar curvature).

In \cite{T} a classification theorem for 4-dimensional conformal flat nearly K\"ahler
manifolds is proved and in \cite{TW} is proved a classification theorem for conformal 
flat nearly K\"ahler manifolds. In \cite{K} we have considered analogous problem for
$M \in AK_3$ and ${\rm dim}\, M \ge 6$. We note that if ${\rm dim}\, M = 4$ then
$M \in NK$ is equivalent to $M \in K$ and $M \in QK$ is equivalent to $M \in AK$.
Now we shall prove:

\vspace{0.1in}
{\bf Theorem.} Let $M \in QK_3$ be conformal flat and ${\rm dim}\, M =4$. Then $M$ is
either of constant sectional curvature or is locally a product of $M_1 \times M_2$, where
$M_1$ (resp. $M_2$) is a 2-dimensional K\"ahler manifold of constant sectional
curvature $c$ (resp. $-c$), $c>0$.                           

\vspace{0.1in}
First we obtain:

\vspace{0.1in}
{\bf Lemma.} Let $M \in QK_3$ be conformal flat. Then $\tau$ is a constant and
$$
	2(\nabla_X S)(Y,Z) = S((\nabla_X J)Y,JZ) + S(JY,(\nabla_X J)Z).   \leqno (2)
$$

\vspace{0.1in}
{\bf Proof.} We note, that $M\in QK_3$ implies that $S(X,Y)=S(JX,JY)$. Let
$X,Y \in \X (M)$, $g(X,Y)=g(X,JY)=0$. According to the second Bianchi identity
$$
	(\nabla_{JX}R)(X,Y,JY,X) + (\nabla_X R)(Y,JX,JY,X) + (\nabla_Y R)(JX,X,JY,X)=0.
$$
Hence, using (1)
$$
	(\nabla_Y S)(X,Y)=S(JX,(\nabla_YJ)Y) + S(JY,(\nabla_YJ)X) - S(Y,(\nabla_{JX}J)Y).
$$
Since $M \in QK_3$ this implies that
$$
	(\nabla_Y S)(X,Y) + (\nabla_{JY} S)(X,JY) =0 .    \leqno (3) 
$$
From (3) we find
$$
	\frac 12 X(\tau) = (\nabla_X S)(X,X) + (\nabla_{JX} S)(X,JX) \leqno (4)
$$
where $g(X,X)=1$. On other hand, from the second Bianchi identity and (1) it follows:
$$
	(\nabla_X S)(Y,Z) - (\nabla_Y S)(X,Z) = \frac 1{2(2m-1)}\{ X(\tau)g(Y,Z) - Y(\tau)g(X,Z) \} \leqno (5)
$$
where ${\rm dim} \ M=2m$. In particular
$$
	X(\tau) = 2(2m-1) \{ (\nabla_X S)(JX,JX)-(\nabla_{JX} S)(X,JX) \}  . \leqno(6)
$$
Since $M \in QK_3$, (4) and (6) imply that $\tau = \rm const.$ and
$$
	(\nabla_X S)(X,X) + (\nabla_{JX} S)(X,JX) = 0 .
$$
Hence, using (3) we conclude that (3) holds for all $X,Y \in \X (M)$. From
$\tau = \rm const.$ and (5) it follows that $\nabla S$ is a symmetric tensor.
Now (3) implies that
$$
	(\nabla_X S)(Y,Z) + (\nabla_{X} S)(JY,JZ) =0  
$$
and hence we obtain (2).

\vspace{0.1in}
{\bf Proof of the theorem.} Let $p \in M$ and $\{ e_1,\, e_2,\, Je_1,\, Je_2 \}$ be an
orthonormal basis of $T_p(M)$ such that $Se_i=\lambda_ie_i$, $i=1,2$. Then
$$
	(\nabla_{e_i}S)(e_j,e_j)=(\nabla_{e_j}S)(e_i,e_j)
$$
and (2) imply
$$
	(\lambda_i-\lambda_j)g((Je_i,(\nabla_{e_j}J)e_j)=0 .
$$
We have two cases:

Case 1. $\lambda_1=\lambda_2$, i.e. $M$ is Einsteinian in $p$.

Case 2. $\lambda_1 \ne \lambda_2$. Then 
$$
	g((Je_i,(\nabla_{e_j}J)e_j)=g(e_i,(\nabla_{e_j}J)e_j)=0
$$
for all $i,j=1,2$.

In both cases (2) implies $\nabla S=0$ in $p$. So $\nabla S=0$.

If $M$ is irreducible, it is an Einsteinian manifold and hence, it is
of constant sectional curvature. If $M$ is reducible, it is locally
a product $M_1 \times M_2$, where $M_1$ and $M_2$ have the stated 
properties.

\vspace {0.7cm}
\begin{flushright}
{\it Faculty of Mathematics and Mechanics \\
University of Sofia \\
Sofia, Bulgaria}
\end{flushright}

\end{document}